\documentclass[preprint,12pt]{elsarticle}
 \journal{}
\usepackage {multicol}
\usepackage{amssymb}
\usepackage[fleqn]{amsmath}
\usepackage{mathrsfs}
\usepackage{hyperref}
\usepackage{graphicx}
\usepackage{color}
\usepackage{caption}
\usepackage{esint}
\usepackage{dsfont}

\begin{document}
\begin{frontmatter}

\title{Extension of Donsker's Invariance Principle with Incomplete Partial-Sum Process}

\author{Jingwei Liu \corref{cor1}}
\ead{liujingwei03@tsignhua.org.cn}
\cortext[cor1]{Corresponding author.}
\address{School of Mathematical Sciences,Beihang University,Beijing,100083,P.R China}

\begin{abstract}
Based on deleting-item central limit theory, the classical Donsker's theorem of partial-sum process of independent and identically distributed (i.i.d.) random variables is extended to incomplete partial-sum process. The incomplete partial-sum process Donsker's invariance principles are constructed and derived for
general partial-sum process of i.i.d random variables and empirical process respectively, they are
not only the extension of functional central limit theory, but also the extension of deleting-item central limit theory. Our work enriches the random elements structure of weak convergence.
\end{abstract}


\end{frontmatter}

\section{Introduction}
\label{}

 Donsker's theorem or Donsker's invariance principle is a well-known functional extension of the central limit theorem in probability theory [1,2]. The original Donsker's invariance principle deals with a sequence of independent and identically distributed (i.i.d.) random variables[1], then the invariance principle is extended to dependent variables,stationary process,uniform mixing,martingale, markov chain, etc. [1-9].
 Donsker's invariance principle for empirical process, which can be expressed as a sum of i.i.d random variables, is an important theorem in probability and statistics[3,9-13].
 The above two Donsker's invariance principles are developed on Skorohod space $D[0,1]$,and $D[0,\infty)$ respectively [2,3,12,13-15], where $D[0,1]$ is the space of functions on [0,1] which are right--continuous and have left--hand limits (functions with jump discontinuities) , and $D[0,\infty)$ is defined in the same way.

The motivation of this paper aims to extend the original Donsker's invariance principle with partial-sum process to Donsker's invariance principle with incomplete partial-sum process. It is the continuous work of deleting-item CLT in our previous work [16].

The rest of the paper are organized as follows. In Section 2, the preliminary notation and classical Donsker's invariance principle are given.
Section 3 constructs and proves the incomplete partial-sum process Donsker's invariance principles. The discussion is given in Section 4.

\section{Notation and Classical Donsker's Theorems}
\label{}

\textbf{Definition 1} Let $\xi_1,\xi_2,\ldots,\xi_n,\ldots$ be a sequence of independent and identically distributed (i.i.d.) random variables, $E\xi_n=0$ and $D\xi_n=\sigma^2<+\infty$. Let
$S_n=\displaystyle \sum\limits_{i=1}^{n}\xi_i$, And,
\begin{equation}
\begin{array}{ll}
 W_n(t)=\displaystyle \frac{1}{\sigma\sqrt{n}} \sum\limits_{i=1}^{[nt]}\xi_i=\displaystyle \frac{1}{\sigma\sqrt{n}} S_{[nt]}, \ t\in [0,1].
\end{array}
\end{equation}
$W_n(t)$ is called partial-sum process ([12],P225).

The classical Donsker invariance principles about general i.i.d random variables and empirical process are as follows[1,2,3,9,12].

\vspace{0.5cm}

\textbf{Theorem 1} Let $\xi_1,\xi_2,\ldots,\xi_n,\ldots$  be a sequence of independent and identically distributed (i.i.d.) random variables,$E\xi_n=0$ and $D\xi_n=\sigma^2<+\infty$. Let
\begin{equation}
\begin{array}{ll}
 \displaystyle W_n(t)&=\displaystyle \frac{1}{\sigma \sqrt{n}}S_{\lfloor nt\rfloor},\ \ t\in [0,1]. \\
 \displaystyle X_n(t)&=\displaystyle \frac{1}{\sigma \sqrt{n}} [S_{[nt]} +(nt-[nt])\xi_{[nt]+1}], \ \ t\in [0,1].
\end{array}
\end{equation}
Let $W=(W(t))_{t\in [0,1]}$ be a standard Brownian motion (or Wiener process). If $W_n(t)$ and $X_n(t)$ are in the Skorokhod space $\mathcal{D}[0,1]$. Then,
\begin{equation}
\begin{array}{ll}
(1)\ \displaystyle W_n(t)\stackrel{d}{\longrightarrow} W(t),\ \ n\to \infty, t\in [0,1].\\
(2)\ \displaystyle X_n(t)\stackrel{d}{\longrightarrow} W(t),\ \ n\to \infty, t\in [0,1].\\
\end{array}
\end{equation}
where ``$\stackrel{d}{\longrightarrow}$'' denotes the convergence in distribution, $[\cdot]$ denotes $[x]=\sup\{m:m\in \mathbb{Z}, m\leq x \},x\in \mathbb{R}$. $X_n(t)$ is the polygonal function [2].
Equivalently ([2], P16, Th2.1), \\
(i) for every bounded and continuous real-valued functional $f$ on the space $C[0,1]$ of continuous functions on the interval [0,1], with the uniform topology, the weakly convergence
\begin{equation}
\begin{array}{ll}
\displaystyle Ef(W_n(t))\rightarrow Ef(W(t)),\ n\rightarrow +\infty,\Longleftrightarrow \displaystyle W_n(t)\stackrel{d}{\longrightarrow} W(t) \\
\displaystyle Ef(X_n(t))\rightarrow Ef(W(t)),\ n\rightarrow +\infty,\Longleftrightarrow \displaystyle X_n(t)\stackrel{d}{\longrightarrow} W(t)
\end{array}
\end{equation}
(ii) for an arbitrary set $G$ in the Borel $\sigma-$ algebra $B_c$ in $C_{[0,1]}$ with
$P\{W\in \partial G\}=0$ (P--continuity set), one has
\begin{equation}
\begin{array}{ll}
\displaystyle P\{W_n(t)\in G\}\rightarrow P\{W(t)\in G\},\ n\rightarrow +\infty, \Longleftrightarrow \displaystyle W_n(t)\stackrel{d}{\longrightarrow} W(t)\\
\displaystyle P\{X_n(t)\in G\}\rightarrow P\{W(t)\in G\},\ n\rightarrow +\infty, \Longleftrightarrow \displaystyle X_n(t)\stackrel{d}{\longrightarrow} W(t)
\end{array}
\end{equation}

\vspace{0.5cm}

\textbf{Definition 2} Let  $\xi_1$, $\xi_2$, $\ldots$,$\xi_n$,$\ldots$ be the sequence of i.i.d. random variables with distribution function $F(x)$. Let
\begin{equation}
\begin{array}{ll}
\displaystyle F_{n}(x)=\displaystyle \frac{1}{n} \sum_{i=1}^n 1_{\{\xi_i\leq x\}}
\end{array}
\end{equation}
where $1_{\{A\}} $ is the indicator function of set $A$. $F_n(x)$ is called empirical distribution function of $\{\xi_n\}$. Define the centered and scaled version of $F_n(x)$ by
\begin{equation}
\begin{array}{ll}
\displaystyle G_{n}(x)={\sqrt {n}}(F_{n}(x)-F(x))
\end{array}
\end{equation}
which is indexed by $x\in R$. $G_{n}(x)$ is called empirical process.

According to the classical central limit theorem, for fixed $x$, the random variable $G_n(x)$ converges in distribution to a Gaussian (normal) random variable $G(x)$ with zero mean and variance
$F(x)(1-F(x))$ as $n\rightarrow +\infty$. The classical Donsker's theorem of empirical process can be expressed as follows[2,9-15].

\vspace{0.5cm}

\textbf{Theorem 2} The sequence of $G_n(x)$, as random elements of the Skorokhod space $\mathcal{D}[0,1]$ (or $\mathcal{D}[0,\infty)$ ), converges in distribution to a Gaussian process $G$ with zero mean and covariance given by
\begin{equation}
\begin{array}{ll}
\displaystyle \operatorname {cov} [G(s),G(t)]=E[G(s)G(t)]=\min\{F(s),F(t)\}-F(s)F(t).
\end{array}
\end{equation}
The process $G(x)$ can be written as $B(F(x))$ where $B$ is a standard Brownian bridge on [0,1] ( or $[0,\infty)$ ).

\vspace{0.5cm}

Generally, let $\xi_1$, $\xi_2$, $\ldots$,$\xi_n$,$\ldots$ be i.i.d. random elements with law $P$ in the measurable space $(\mathcal{X},\mathcal{A})$, and let $\mathcal{F}$ be a collection of square-integrable, measurable functions $f:\mathcal{X}\mapsto \mathbb{R}$. The sequential empirical process [12] is defined as
\begin{equation}
\begin{array}{ll}
\displaystyle \mathbb{Z}_{n}(x,f)=\displaystyle \frac{1}{\sqrt{n}}\sum\limits_{i=1}^{[nx]} (f(\xi_i)-Pf)=\displaystyle \sqrt{\frac{[nx]}{n}} \mathbb{G}_{[nx]}(f). \\
\end{array}
\end{equation}
where $\displaystyle G_{n}(x)={\sqrt {n}}(P_{n}-P)$ is the empirical process indexed by $\mathcal{F}$. The index $(s,f)$ ranges over $[0,1]\times \mathcal{F}$. And, $Pf=\int_{\mathcal{X}}fdP$.

\vspace{0.5cm}

\textbf{Theorem 3}  The marginals of sequence of processes $\{\mathbb{Z}_{n}(x,f):(x,f)\in [0,1]\times \mathcal{F}\}$ converge to the marginals of a Gaussian process $\mathbb{Z}$ (known as Keifer--M$\ddot{u}$ller process), with zero mean and covariance given by
\begin{equation}
\begin{array}{ll}
\displaystyle \operatorname {cov} (\mathbb{Z}(s,f),\mathbb{Z}(t,g))=(s\wedge t)(Pfg-PfPg).
\end{array}
\end{equation}

\section{Extension of Donsker's Invariance Principles}

Suppose $\xi_1,\xi_2,\cdots,\xi_n,\cdots$ is a real-valued random variable sequence in $\mathbf{R}$. We only address $\mathbf{R}$ case in this paper.
Denote $J_n=\{1,2,\cdots,n\}$, $\{i_1,i_2,\ldots,i_{n}\}$ is an arbitrary permutation of $J_n$.
Denote $J_{k^*}=\{i_1,i_2,\cdots,i_{k^*}\}$, where $\{i_1,i_2,\cdots,i_{k^*}\}$ is any $k^*$ different elements of $J_n$. Obviously, there are $C_n^{k^*}$ possible combination of $J_{k^*}$.
Denote\\
\begin{equation}
\begin{array}{ll}
S_n=S_{J_n}=\sum\limits_{i=1}^n \xi_i, \\
\tilde{S}_{(n,{k^*})}=S_{J_n\backslash J_{k^*}}=\sum\limits_{i\in J_n \backslash J_{k^*}} \xi_i .
\end{array}
\end{equation}
For convenience, Denote $J_n=\{1,2,\cdots,n\}$, and
\begin{equation}
k^{*}\stackrel{\bigtriangleup}{=} k^{*}_{[n]}=\left\{
\begin{array}{ll}
k,    & 0\leq k<n. \\
k(n), & 0\leq k(n)< n.\\
\end{array}
\right.
\  \ J_{k^{*}}=\left\{
\begin{array}{ll}
J_{k},   & k^{*}=k   \\
J_{k(n)},& k^{*}=k(n). \\
\end{array}
\right.
\end{equation}

If
\begin{equation}
\begin{array}{ll}
\displaystyle \lim \limits_{n\rightarrow \infty} \frac{k^{*}}{n}=0,
\end{array}
\end{equation}
it is called asymptotic deleting negligibility condition [16]. Note that $k^{*}$ could be infinite as $n\rightarrow +\infty$, for example, $k^{*}=[n^r]$, ($0<r<1$).

The index of $(n,k^{*})$ is adopted here for deleting-item meaning is different from the $(n,k(n))$ for traditional triangular array in probability theory [2].
$S_{(n,k^*)}$ is somewhat like a constrained 2-dimension lattice stochastic process, they are generated by  considering all possible combination of $C_n^{k^*}$ for each $n$. As $n\rightarrow +\infty$, the possible combination of $S_{(n,k^*)}$ will be countable infinite. Hence, we obtain the following conclusion.

\textbf{Lemma 1} Let $\xi_1$, $\xi_2$, $\ldots$,$\xi_n$,$\ldots$ be an i.i.d. random variable sequence, $E\xi_k=0$, $D\xi_k=\sigma^2$. Then, $S_{(n,k^*)}$ will be generated countable infinite by $\{\xi_n\}$ as
$n\rightarrow +\infty$, constrained by $\displaystyle \lim\limits_{n\rightarrow \infty} \frac{k^{*}}{n}=0$.

\vspace{0.5cm}

\textbf{Definition 3} Let $\xi_1,\xi_2,\ldots,\xi_n,\ldots$ be a sequence of independent and identically distributed (i.i.d.) random variables, $E\xi_n=0$ and $D\xi_n=\sigma^2<+\infty$. Suppose $\{i_1,i_2,\ldots,i_{[nt]}\}$ is an arbitrary permutation of $J_{[nt]}=\{1,2,\cdots,[nt]\}$, for a natural number $k^{*}$, $(0<k^{*}<[nt])$, denote $J_{k^{*}_{[nt]}}=\{i_1,i_2,\ldots,i_k^{*}\}$,
\begin{equation}
\begin{array}{ll}
\displaystyle \tilde{S}_{(n,k^{*})}(t)=\displaystyle S_{{J_{[nt]}\backslash J_{k^{*}_{[nt]}}}}
=\displaystyle \sum\limits_{i\in J_{[nt]}\setminus J_{k^{*}_{[nt]}}} \xi_i, \  t\in [0,1].\\
\displaystyle \tilde{W}_{(n,k^{*})}(t)=\displaystyle \frac{1}{\sigma \sqrt{n}} \sum\limits_{i\in J_{[nt]}\setminus J_{k^{*}_{[nt]}}} \xi_i=\displaystyle \frac{1}{\sigma \sqrt{n}}\tilde{S}_{(n,k^{*})}(t), \  t\in [0,1].\\
\end{array}
\end{equation}
$\tilde{W}_{(n,k^{*})}(t)$ is called deleting-item partial-sum process, it is also an incomplete partial-sum process. Denote
\begin{equation}
\begin{array}{ll}
\displaystyle \tilde{X}_{(n,k^{*})}(t)=\displaystyle \frac{1}{\sigma \sqrt{n}} \{S_{{J_{[nt]}\backslash J_{k^*_{[nt]}}}} +(nt-[nt])\xi_{[nt]+1} \},  (0\leq t\leq1),
\end{array}
\end{equation}
which is the deleting-item polygonal function [2].

\vspace{0.5cm}

Since the original idea of invariance principle first realized by A. Kolmogorov (1931) is to compute the limit distribution by using a special case and then passing to the general case, without a loss of generality, to distinguish the partial-sum process and incomplete partial-sum process apart, Lemma 1 investigates their stochastic property in discrete random variables case.

\textbf{Lemma 2} Let $\xi_1$, $\xi_2$, $\ldots$,$\xi_n$,$\ldots$ be discrete i.i.d. random variables, $E\xi_k=0$, $D\xi_k=\sigma^2$, $T\subseteq R$ ($T=[0,1]$ or $[0,+\infty)$ ). Then,
\begin{itemize}
\item[(1)]$\{S_n(t)\}_{t\in T}$ is an incremental process.
\item[(2)]$\{S_{J\setminus J_{k^*}}\}_{t\in T}$ are generally not independent incremental processes.
\item[(3)]$\{\tilde{S}_{(n,k^*)}(t)\}_{t\in T}$ is not stationary, not Markov process, not ergodic process.
\item[(4)]$\{\tilde{S}_{(n,k^*)}(t)\}_{t\in T}$ are infinite stochastic processes different from $\{S_n(t)\}_{t\in T}$.
\item[(5)]$\{\tilde{W}_{(n,k^*)}(t)\}_{t\in T}$ are infinite stochastic processes different from $\{W_n(t)\}_{t\in T}$.
\item[(6)]$\{\tilde{X}_{(n,k^*)}(t)\}_{t\in T}$ are infinite stochastic processes different from $\{X_n(t)\}_{t\in T}$.
\end{itemize}

\textbf{Proof}

(1) This conclusion is obvious. For $\forall t_0,t_1,\cdots,t_m\in T$,satisfying $t_0<t_1<\cdots<t_m$,
\begin{equation}
\begin{array}{ll}
\{S_n(t_0), S_n(t_1)-S_n(t_0), \cdots , S_n(t_{m})-S_n(t_{m-1})\} \\
=\{\displaystyle \sum\limits_{i\in J_{[nt_0]}} \xi_i, \sum\limits_{i\in J_{[nt_1]}\setminus J_{[nt_0]}} \xi_i, \cdots , \sum\limits_{i\in J_{[nt_m]}\setminus J_{[nt_{m-1}]}} \xi_i \}
\end{array}
\end{equation}
Since $\{\xi_n\}$ is i.i.d sequence, $S_n(t_0),S_n(t_1)-S_n(t_0),\cdots,S_n(t_{m})-S_n(t_{m-1})$ are independent. Hence, $\{S_n(t)\}_{t\in T}$ is an independent incremental process.

(2)For $\forall t_0,t_1,\cdots,t_m\in T$, satisfying $t_0<t_1<\cdots<t_m$,
\begin{equation}
\begin{array}{ll}
\displaystyle \{\tilde{S}_{(n,k^*)}(t_0), \tilde{S}_{(n,k^*)}(t_1)-\tilde{S}_{(n,k^*)}(t_0), \cdots , \tilde{S}_{(n,k^*)}(t_{m})-\tilde{S}_{(n,k^*)}(t_{m-1})\} \\
=\{\displaystyle \sum\limits_{i\in J_{[nt_0]}\setminus J_{k^*_{[nt_0]}}} \xi_i,
 \sum\limits_{i\in \{J_{[nt_1]}\setminus J_{k^*_{[nt_1]}}\} \setminus \{J_{[nt_0]}\setminus J_{k^*_{[nt_0]}}\} } \xi_i, \cdots ,\\
\ \ \ \ \sum\limits_{i\in \{J_{[nt_m]}\setminus J_{k^*_{[nt_m]}}\} \setminus \{J_{[nt_{m-1}]} \setminus J_{k^*_{[nt_{m-1}]}}\}} \xi_i \}
\end{array}
\end{equation}
Denote\\
\begin{equation}
\begin{array}{ll}
 I_0= \{J_{[nt_0]}\setminus J_{k^*_{[nt_0]}}\},\\
 I_i= \{J_{[nt_i]}\setminus J_{k^*_{[nt_i]}}\} \setminus \{J_{[nt_{i-1}]}\setminus J_{k^*_{[nt_{i-1}]}}\}.\ i=1,2,\cdots,m.
\end{array}
\end{equation}
As in $\displaystyle \tilde{S}_{(n,k^*)}(t)$, $k^*=k^*_{[nt]}$, and the choice of $J_{k^*_{[nt]}}$ is arbitrary in $J_{[nt]}$, which leads to $I_i \bigcap I_{i-1} \neq \emptyset$.
Then $\operatorname {cov}(\sum\limits_{j\in I_{i-1}}\xi_j,\sum\limits_{j\in I_{i}}\xi_j)\neq 0$.
Therefore, the conclusion holds.

(3) Since
$ \displaystyle D\tilde{S}_{(n,k^*)}(t)=\displaystyle (1-\frac{k^*_{[nt]}}{[nt]})\sigma^2 $
is determined by the choice of $k^*_{[nt]}$, $\tilde{S}_{(n,k^*)}(t)$ is not ergodic process.
Again, $\displaystyle \operatorname {cov}(\tilde{S}_{(n,k^*)}(t), \tilde{S}_{(n,k^*)}(s))$
is affected by $J_{k^*_{[nt]}}$ and $J_{k^*_{[ns]}}$ and the choice of $k^*_{[nt]}$ and $k^*_{[ns]}$, $\tilde{S}_{(n,k^*)}(t)$ is not stationary process.

For $\forall t_0<t_1<\cdots,<t_m<s \in T$, and possible status $i_0,i_1,\cdots,i_m,i$,
for $s<t+s\in T$, since the deleting-item way of $\tilde{S}_{(n,k^*)}(\cdot)$,

\begin{equation*}
\begin{array}{ll}
\displaystyle P(\tilde{S}_{(n,k^*)}(t+s)=j|\tilde{S}_{(n,k^*)}(t_0)=i_0,\tilde{S}_{(n,k^*)}(t_1)=i_1,\cdots,\tilde{S}_{(n,k^*)}(t_m)=i_m,\\
\ \ \ \ \tilde{S}_{(n,k^*)}(s)=i) \neq P(\tilde{S}_{(n,k^*)}(t+s)=j|\tilde{S}_{(n,k^*)}(s)=i)
\end{array}
\end{equation*}
$\tilde{S}_{(n,k^*)}(t)$ is not a Markov chain.

(4) For fixed $t>0$ and each $n$ large enough, there are $C_{[nt]}^{k^*_{[nt]}}$ possible combinations (deleting-item style), even one combination way is chosen, the total ways to determine a deleting-item partial-sum process will be at least $\displaystyle \prod\limits_{i=n}^{+\infty} C_{[it]}^{k_{[it]}^*}$. Hence, from the original i.i.d random variable sequence, infinite deleting-item partial-sum processes are generated different from partial-sum process.

(5) and (6) hold according to (4).

Hence, it ends the proof.

\hfill $\Box$

\textbf{Lemma 3} $\forall s<t\in T$,
\begin{equation}
\begin{array}{ll}
(1) \displaystyle (\tilde{W}_{(n,k^*)}(s),\tilde{W}_{(n,k^*)}(t)-\tilde{W}_{(n,k^*)}(s)) \stackrel{d}{\longrightarrow} (W(s),W(t)-W(s)).\\
(2) \displaystyle (\tilde{X}_{(n,k^*)}(s),\tilde{X}_{(n,k^*)}(t)-\tilde{X}_{(n,k^*)}(s)) \stackrel{d}{\longrightarrow} (W(s),W(t)-W(s)).\\
\end{array}
\end{equation}

\textbf{Proof}

(1) According to the conclusion in Donsker's theorem [2], Chebyshev inequality and Slutsky's theorem,
\begin{equation}
\begin{array}{ll}
\displaystyle  (\tilde{W}_{(n,k^*)}(s),\tilde{W}_{(n,k^*)}(t)-\tilde{W}_{(n,k^*)}(s))\\
=\displaystyle  (\frac{1}{\sigma \sqrt{n}}\tilde{S}_{(n,k^*)}(s), \frac{1}{\sigma \sqrt{n}}\tilde{S}_{(n,k^*)}(t)-\frac{1}{\sigma \sqrt{n}}\tilde{S}_{(n,k^*)}(s) )\\
=\displaystyle  (\frac{1}{\sigma \sqrt{n}}S_{[ns]},
\frac{1}{\sigma \sqrt{n}}S_{[nt]}-\frac{1}{\sigma \sqrt{n}}S_{[ns]})\\
\ \ \ \ \displaystyle  -(\frac{1}{\sigma \sqrt{n}} \sum\limits_{i\in J_{k^*_{[ns]}}} \xi_i,
\frac{1}{\sigma \sqrt{n}} \sum\limits_{i\in J_{k^*_{[nt]}}} \xi_i
-\frac{1}{\sigma \sqrt{n}} \sum\limits_{i\in J_{k^*_{[ns]}}} \xi_i)\\
\displaystyle
\stackrel{d}{\longrightarrow} (W(s),W(t)-W(s))+(0,0) \\
\stackrel{d}{\longrightarrow} (W(s),W(t)-W(s)).
\end{array}
\end{equation}

(2) Similar to the above proof, the conclusion holds.

\hfill $\Box$

Lemma 2 shows that, though $\displaystyle \tilde{W}_{(n,k^*)}(t)$ and $\displaystyle \tilde{X}_{(n,k^*)}(t)$ are not independent incremental process, their limit stochastic processes are still standard Brownian motion, which has independent increment property.

\vspace{0.5cm}

The deleting-item Donsker's invariance principles are given below .

\textbf{Theorem 4} Let $\xi_1,\xi_2,\ldots\xi_n,\ldots$ be a sequence of independent and identically distributed (i.i.d.) random variables with $E\xi_n=0$ and $D\xi_n=\sigma^2$. If
$\displaystyle \lim \limits_{n\rightarrow \infty} \frac{k^{*}}{n}=0$, then
\begin{equation}
\begin{array}{ll}
(1)\ \displaystyle \tilde{W}_{(n,k)}(t)\stackrel{d}{\longrightarrow} W(t),\ \ n\to \infty, t\in [0,1].\\
(2)\ \displaystyle \tilde{X}_{(n,k)}(t)\stackrel{d}{\longrightarrow} W(t),\ \ n\to \infty, t\in [0,1].\\
\end{array}
\end{equation}
Equivalently, \\
(i) for every bounded and continuous real-valued functional $f$ on the space $C[0,1]$ of continuous functions on the interval [0,1], with the uniform topology, the weakly convergence
\begin{equation}
\begin{array}{ll}
\displaystyle Ef(\tilde{W}_{(n,k)}(t))\rightarrow Ef(W(t)),\ n\rightarrow +\infty,\Longleftrightarrow \displaystyle \tilde{W}_{(n,k)}(t)\stackrel{d}{\longrightarrow} W(t) \\
\displaystyle Ef(\tilde{X}_{(n,k)}(t))\rightarrow Ef(W(t)),\ n\rightarrow +\infty,\Longleftrightarrow \displaystyle \tilde{X}_{(n,k)}(t)\stackrel{d}{\longrightarrow} W(t)
\end{array}
\end{equation}
(ii) for an arbitrary set $G$ in the Borel $\sigma-$ algebra $B_c$ in $C_{[0,1]}$ with
$P\{W\in \partial G\}=0$, one has
\begin{equation}
\begin{array}{ll}
\displaystyle P\{\tilde{W}_{(n,k)}(t)\in G\}\rightarrow P\{W(t)\in G\},\ n\rightarrow +\infty, \Longleftrightarrow \displaystyle \tilde{W}_{(n,k)}(t)\stackrel{d}{\longrightarrow} W(t)\\
\displaystyle P\{\tilde{X}_{(n,k)}(t)\in G\}\rightarrow P\{W(t)\in G\},\ n\rightarrow +\infty, \Longleftrightarrow \displaystyle \tilde{X}_{(n,k)}(t)\stackrel{d}{\longrightarrow} W(t)
\end{array}
\end{equation}

\textbf{Proof}.

(1) According to Theorem 1(2) and deleting-item central limit theorem [16],
$\displaystyle \tilde{W}_{(n,k)}(t)\stackrel{d}{\longrightarrow} \displaystyle W(t)$.

(2) According to Theorem 1(3) and deleting-item central limit theorem [16],
$\displaystyle \tilde{X}_{(n,k)}(t)\stackrel{d}{\longrightarrow} \displaystyle W(t)$.

And the following conclusions about $\displaystyle Ef(\tilde{W}_{(n,k)}(t))$,
$\displaystyle Ef(\tilde{X}_{(n,k)}(t))$, \\
$\displaystyle P\{ \tilde{W}_{(n,k)} (t) \in G\}$ and $\displaystyle P\{\tilde{X}_{(n,k)} (t)\in G\}$ are obvious based on the theory of weak convergence and convergence in distribution ([2],P16,Th2.1).

\hfill $\Box$.

\textbf{Definition 4} Let $\xi_1$, $\xi_2$, $\ldots$,$\xi_n$,$\ldots$ be i.i.d. random variables with distribution function $F(x)$, Define
\begin{equation}
\begin{array}{ll}
\displaystyle \hat{F}_{(n,k^*)}(x)=\displaystyle \frac{1}{n} \sum\limits_{i \in J_{n}\setminus J_{k^*_{n}} } 1_{\{\xi_i\leq x\}}\\
\displaystyle \tilde{F}_{(n,k^*)}(x)=\displaystyle \frac{1}{n} \sum\limits_{i \in J_{n}\setminus J_{k^*_{n}} } (1_{\{\xi_i\leq x\}}-F(x))\\
\end{array}
\end{equation}
where $1_{\{A\}} $ is the indicator function of set $A$, $1\leq k^*_{n} \leq n$.  $\hat{F}_{(n,k^*)}(x)$ is called deleting-item (or incomplete) empirical distribution function, $\tilde{F}_{(n,k^*)}(x)$ is called deleting-item (or incomplete) centered  empirical distribution function .
Define the deleting-item centered and scaled version of $\tilde{F}_{(n,k^*)}(x)$ by
\begin{equation}
\begin{array}{ll}
\displaystyle \tilde{G}_{(n,k^*)}(x)={\sqrt {n}}(\tilde{F}_{(n,k^*)}(x))
\end{array}
\end{equation}
indexed by $x\in R$. $\tilde{G}_{(n,k^*)}(x)$ is called deleting-item empirical process.\footnote{The definition of $\tilde{G}_{(n,k^*)}(x)$ here is to keep uniform formula for [nx] case especially in deleting-item Keifer-M$\ddot{u}ller$ process.}

\textbf{Lemma 4} Constrained by $\displaystyle \lim\limits_{n\rightarrow \infty} \frac{k^{*}}{n}=0$, one has \\ (1)$\tilde{F}_{(n,k^*)}(x)$ are generated infinite stochastic processes indexed by $(n,k^*)$.\\
(2)$\tilde{G}_{(n,k^*)}(x)$ are generated infinite stochastic processes indexed by $(n,k^*)$.

\textbf{Theorem 5.} If
$\displaystyle \lim \limits_{n\rightarrow \infty} \frac{k^{*}}{n}=0$, then the sequence of $\tilde{G}_{(n,k^*)}(x)$, as random elements of the Skorokhod space $\mathcal{D}[0,1]$ (or $\mathcal{D}[0,\infty)$ ), converges in distribution to a Gaussian process $G$ with zero mean and covariance given by
\begin{equation}
\begin{array}{ll}
\displaystyle \operatorname {cov} [G(s),G(t)]=E[G(s)G(t)]=\min\{F(s),F(t)\}-F(s)F(t).
\end{array}
\end{equation}
The process $G(x)$ can be written as $B(F(x))$ where $B$ is a standard Brownian bridge on [0,1] (or $[0,\infty)$).

\textbf{Proof}
According to Theorem 2,
\begin{equation}
\begin{array}{ll}
\displaystyle G_{n}(x)&=\displaystyle {\sqrt {n}}(F_{n}(x)-F(x))=\displaystyle{\sqrt {n}}(\frac{1}{n} \sum_{i=1}^n 1_{\{\xi_i\leq x\}}-F(x))\\
&=\displaystyle \frac{1}{\sqrt {n}} \sum_{i=1}^n (1_{\{\xi_i\leq x\}}-F(x))\\
                      &\stackrel{d}{\longrightarrow} B(F(x)),\ \ n\to \infty. \\
\end{array}
\end{equation}
According to Slutsky's Theorem or deleting-item classical central limit theorem [16],
\begin{equation}
\begin{array}{ll}
\displaystyle \tilde{G}_{(n,k^*)}(x)&=\displaystyle {\sqrt {n}}\tilde{F}_{(n,k^*)}(x)\\
                                    &=\displaystyle {\sqrt {n}}\displaystyle \frac{1}{n} \sum\limits_{i \in J_{n}\setminus J_{k^*_{n}} } (1_{\{\xi_i\leq x\}}-F(x))\\
&\stackrel{d}{\longrightarrow} B(F(x)),\ \ n\to \infty.
\end{array}
\end{equation}
Hence, it ends the proof.

\hfill $\Box$.

\vspace{0.5cm}

\textbf{Definition 5} let $\xi_1$, $\xi_2$, $\ldots$,$\xi_n$,$\ldots$ be i.i.d. random elements with law $P$ in the measurable space $(\mathcal{X},\mathcal{A})$, and let $\mathcal{F}$ be a collection of square-integrable, measurable functions $f:\mathcal{X}\mapsto \mathbb{R}$.
The deleting-item sequential empirical process is defined as
\begin{equation}
\begin{array}{ll}
\displaystyle \tilde{\mathbb{Z}}_{(n,k^*)}(x,f)=\displaystyle \frac{1}{\sqrt{n}}\sum\limits_{i \in J_{[nx]}\setminus J_{k^*_{[nx]}} } (f(\xi_i)-Pf)=\displaystyle \sqrt{\frac{[nx]}{n}} \tilde{\mathbb{G}}_{(n,k^*)}(f).
\end{array}
\end{equation}
where $\displaystyle \tilde{\mathbb{G}}_{(n,k^*)}={\sqrt {n}}(\tilde{P}_{(n,k^*)})$ is the deleting-item empirical process indexed by $\mathcal{F}$, $\displaystyle \tilde{P}_{(n,k^*)}=\displaystyle \frac{1}{n} \sum\limits_{i \in J_{[nx]}\setminus J_{k^*_{[nx]}} } (\delta_{\xi_i}-P)$, $\displaystyle\delta_{\xi_i}$ is the dirac measure. The index $(s,f)$ ranges over $[0,1]\times \mathcal{F}$.

\vspace{0.5cm}

\textbf{Theorem 6}  If
$\displaystyle \lim \limits_{n\rightarrow \infty} \frac{k^{*}}{n}=0$, then  the marginals of sequence of processes $\{\tilde{\mathbb{Z}}_{(n,k^*)}(x,f)\in [0,1]\times \mathcal{F}\}$ converge to the marginals of a Gaussian process $\mathbb{Z}$ (known as Keifer--M$\ddot{u}$ller process), with zero mean and covariance given by
\begin{equation}
\begin{array}{ll}
\displaystyle \operatorname {cov} (\mathbb{Z}(s,f),\mathbb{Z}(t,g))=(s\wedge t)(Pfg-PfPg).
\end{array}
\end{equation}

\textbf{Proof}
According to Theorem 3,

\begin{equation}
\begin{array}{ll}
\displaystyle \mathbb{Z}_{n}(x,f)=\displaystyle \frac{1}{\sqrt{n}}\sum\limits_{i=1}^{[nx]} (f(\xi_i)-Pf)\stackrel{d}{\longrightarrow} \mathbb{Z}(x,f)

\end{array}
\end{equation}

According to deleting-item classical central limit theorem [16] and
$\displaystyle \lim \limits_{n\rightarrow \infty} \frac{k^{*}}{n}=0$,

$\tilde{\mathbb{Z}}_{(n,k^*)}(x,f)=\displaystyle \frac{1}{\sqrt{n}}\sum\limits_{i \in J_{[nx]}\setminus J_{k^*_{[nx]}} } (f(\xi_i)-Pf) \stackrel{d}{\longrightarrow} \mathbb{Z}(x,f) $.

According to Theorem 3, Keifer--M$\ddot{u}$ller process $\mathbb{Z}$ with zero mean and covariance given by

\begin{equation}
\begin{array}{ll}
\displaystyle \operatorname {cov} (\mathbb{Z}(s,f),\mathbb{Z}(t,g))=(s\wedge t)(Pfg-PfPg).
\end{array}
\end{equation}
which ends the proof.

\hfill $\Box$.

\section{Conclusion}
\label{}
Based on deleting-item central limit theory of i.i.d random variables, infinite countable incomplete partial-sum processes are constructed different from the classical partial-sum process, which converge to Donsker's limit. Hence, enlarging the random elements structure in the previous skorokhod space. Our research reveals the deep complexity of stochastic process structure, and elevates the difficulty in stochastic simulation, especially in stochastic finance. The future work will focus on the Donsker's invariance principles of deleting-item (incomplete) partial-sum process of dependent random variable sequence and application of Donsker's invariance principle in stochastic finance.

\section*{References}
\label{}

\end{document}